\documentclass{article}

\usepackage{amssymb}
\usepackage[cp1251]{inputenc}
\usepackage[english]{babel}
\usepackage{amsmath}
\usepackage{color}

\def\mapr#1{\smash{\mathop{\buildrel{#1}\over\longrightarrow}}}

\def\mapd#1{\Big\downarrow\rlap{$\vcenter{\hbox{$#1$}}$}}

\newtheorem{theorem}{Theorem}
\newtheorem{lemma}{Lemma}
\newtheorem{definition}{Definition}

\newtheorem{cor}{Corollary}
\def\qed{\hfill\vrule width2mm height2mm depth2mm}
\def\proof{{\bf Proof.}}%\nobreak\noindent}
\def\f#1{{$#1$}}
\def\ff#1{{$$#1$$}\noindent}
\def\feq#1#2{\begin{equation}\label{#1}
    #2\end{equation}\noindent}
\def\({\left(}
\def\){\right)}

\def\U{{\bf U}}

\def\id{{\hbox{\bf Id}}}
\def\1{{\bf 1}}

\title{Description of  $G$-bundles over
$G$-spaces with quasi-free proper action of discrete group II}

\author{Morales Mel\'endez\thanks{Partially supported
by the grant 299388 of the mexican 
National Council for Science and
Technology (CONACyT) }, Quitzeh}
\begin{document}
\maketitle

\section{The setting of the problem}

This problem naturally arises from the Connor-Floyd's description
of the bordisms with the action of a group
$G$ using the so-called fix-point construction
This construction reduces the problem of describing the bordisms to two simpler
problems: a) description of the fixed-point set (or, more generally, the stationary point set),
which happens to be a submanifold attached with the structure of its normal bundle
and the action of the same group $G$, however, this action could have stationary points
of lower rank; b) description of the bordisms of lower rank with an action of the group $G$.
We assume that the group  $G$ is discrete.

Lets  $\xi$ be an $G$-equivariant vector bundle with base $M$.
\ff{
\begin{array}{c}
          \xi \\
          \mapd{} \\
          M \\
        \end{array}
}Where the action of the group $G$ is quasi-free over the base with normal
stationary subgroup $H<G$ and there is no more fixed points
of the action of the group $H$ in the total space of the bundle $\xi$.

According \cite[p.1]{MishQuit} the bundle $\xi$ separates as the sum
of its  $G$-subbundles:
\ff{\xi\approx \bigoplus_{k}\xi_k}
where the index runs over all (unitary) irreducible representations
$\rho_{k}: H\mapr{}\U(V_{k})$ of the group $H$ and, as a $H$-bundle $\xi_k$
can be presented as the tensor product:
\ff{
\xi \approx \bigoplus_{k}\eta_k\bigotimes V_{k},
}where
the action of the group $H$ over the bundles $\eta_k$ is trivial,   $V_{k}$
denotes the trivial bundle with fiber $V_{k}$ and with  fiberwise
action of the group $H$, defined using the linear representation $\rho_{k}$.

The particular case $\xi=\eta_{k}\bigotimes V_{k}$ was described in the previous article
\cite{MishQuit}. According \cite[p.14]{MishQuit} the bundle $\xi_k$ can be obtained as the
inverse image of a mapping
\ff{f_k: M/G_0\longrightarrow B\mathrm{Aut}_{G}\(X_k\)}where
$G_0=G/H,\;X_k=G_0\times \(F_k\otimes V_k\)$ is the canonical model and
\f{\mathrm{Aut}_G\(X_k\)} is the group of equivariant automorphisms of the space
 $X_k$ as a vector $G$-bundle
over the base  $G_0$.
So, the bundle $\xi$ can be given by a mapping
\ff{f: M/G_0\longrightarrow \underset{k}{\prod} B\mathrm{Aut}_{G}\(X_k\).}

Consider the vector bundle over the discrete base $G_0$
\feq{can}{X_\rho=G_0\times\( \bigoplus_k\(F_k\otimes V_k\)\).}
Define a fiberwise action $G\times X_\rho\to X_\rho$ by the formula
\ff{
\phi([g],g_{1}):[g]\times \( \bigoplus_k\(F_k\otimes V_k\)\) \to [g_1g]\times \( \bigoplus_k\(F_k\otimes V_k\)\)
}
\feq{canaction}{
\begin{array}{ll}
\phi([g],g_{1})
=
\bigoplus_k\(\id\otimes\rho_k(u(g_{1}g)u^{-1}(g))\)=
\rho(u(g_{1}g)u^{-1}(g)).
    \end{array}
}
\begin{definition}
The bundle \f{X_\rho\mapr{}G_0} with the just defined action is called the
canonical model for the representation $\rho$.
\end{definition}

By  \f{\mathrm{Aut}_G\(X_\rho\)} we denote the group of equivariant automorphisms
of the canonical model  $X_\rho$ as a vector $G$-bundle
over the base  $G_0$ with fiber $\bigoplus_k\(F_k\otimes V_k\)$ and
canonical action of the group $G$.

\begin{lemma}
There exists a monomorphism
\ff{ i:\mathrm{Aut}_{G}\(X_\rho\)\mapr{}\underset{k}{\prod} \mathrm{Aut}_{G}\(X_k\)}
\end{lemma}
\proof
As before, an element of the group  \f{\mathrm{Aut}_G\(X_\rho\)}
is an equivariant mapping$\mathbf{A}^{a}$ such that the pair $(\mathbf{A}^{a},a)$
defines a commutative diagram
\ff{
\begin{array}{ccc}
   X_\rho &\mapr{\mathbf{A}^{a}}& X_\rho\\
  \mapd{} &&\mapd{}\\
   G_{0}&\mapr{a}&G_{0}, \\
\end{array}
}\
\ff{a\in \mathrm{Aut}_G(G_0)\approx G_0,\quad a[g]=[ga],\;[g]\in G_0.}

By the lemma 1 \cite{MishQuit} applied to
group of automorphisms $\mathrm{Aut}_{G}\(X_\rho\)$,
for  $A^a\in \mathrm{Aut}_{G}\(X_\rho\)$,
we have
\ff{A^a|_{X_k}:X_k\mapr{} X_k.}Note that
this restriction is $G$-equivariant.

Define
\ff{ i:\mathrm{Aut}_{G}\(X_\rho\)\mapr{}\underset{k}{\prod} \mathrm{Aut}_{G}\(X_k\)}
by the formula
\ff{i(A^a)=(A^a|_{X_k})_k.}This is clearly a  homomorphism:
it is a product of restrictions over invariant subspaces.
Lets prove that it is injective.
If $A^a|_{X_k}=\id_{X_k}$,
then $A^a=\id_X$ because
\f{X=\oplus_kX_k}.

In order to prove that the image of $i$ is closed, note that it
coincides with those automorphisms which commute with
the inclusion
\ff{\Delta\times \id:X_\rho\mapr{}\underset{k}{\prod} X_k}i.e.
if an element
$(A_k^{a_k})_k\in\underset{k}{\prod} \mathrm{Aut}_{G}\(X_k\) $
leaves the image of \f{X_\rho} invariant, then its
restriction defines an element in
$\mathrm{Aut}_{G}\(X_\rho\)$ and the diagram
\ff{\begin{array}{ccc}
X_\rho&\mapr{A^a}&X_\rho\\
\mapd{}&         &\mapd{}\\
\underset{k}{\prod} X_k&\mapr{(A_k^{a_k})_k}&\underset{k}{\prod} X_k\\
\end{array}}commutes, i.e.  \f{a^k=a\,\forall k} and
\f{A^a|_{X_k}=A_k^a}. In other words
\ff{\prod_kpr_k\circ i (A^a)=\Delta(a)}where
\ff{pr_k: \mathrm{Aut}_{G}\(X_k\)\mapr{} G_0} is the epimorphism of
lemma 2 \cite[p. 9]{MishQuit} and
\ff{\Delta:G_0\hookrightarrow \underset{k}{\prod} G_0.}Then
\ff{i\(\mathrm{Aut}_{G}\(X_\rho\)\)=\prod_kpr_k^{-1}\Delta(G_0).}
\qed

\begin{cor} It takes place an exact sequence of groups
\ff{1\to \prod_k GL(F_k)
\mapr{\varphi}\mathrm{Aut}_G\(X_\rho\)
\mapr{pr} G_0\to 1}
\end{cor}
\proof
Define  $pr=\Delta^{-1}\prod_kpr_k\circ i$.
This is an epimorphism: let $A_k^a\in pr_k(a)$, then
$(A_k^a)_k\in i\(\mathrm{Aut}_{G}\(X_\rho\)\)$.

We define the monomorphism
 \ff{ \varphi:\prod_k GL(F_k) \mapr{}\mathrm{Aut}_G\(X_\rho\)}using
 the monomorphisms
\ff{\varphi_k:GL(F_k)
\mapr{}\mathrm{Aut}_G\(X_k\)}by the formula
 \ff{ \varphi= i^{-1}\prod_k\varphi_k.}Then
 \f{pr\circ\varphi=1} and
 if $a=1$, then $\prod_kpr_k\circ i(A^1)=\Delta(1)$. This means that, for
 every $k$, there is a $B\in GL(F_k)$ such that
 $\varphi_k(B_k)=A^1|_{X^k}$, i.e.
 $\varphi(B_k)_k=A^1$.\qed

Denote by  $\mathrm{Vect}_G(M,\rho)$ the category of  $G$-equivariant
vector bundles $\xi$ over the base $M$ with quasi-free action of the
group $G$ over the base and normal stationary subgroup $H<G$.

Then, by lemma 1 and the observations on p. 13 in \cite{MishQuit}
in terms of homotopy we have
\feq{vectsum}{\mathrm{Vect}_G(M,\rho)\approx \prod_k[M,BU(F_k)]}

Denote by $\mathrm{Bundle}(X,L)$ the category of principal $L$-bundles over
the base $X$.

\begin{theorem} There exists an inclusion
\feq{monocat}{\mathrm{Vect}_G(M,\rho)\longrightarrow
\mathrm{Bundle}(M/G_0,\mathrm{Aut}_G\(X_\rho\)). }
\end{theorem}

\proof\, We already have a monomorphism
\ff{\mathrm{Vect}_G(M,\rho)\longrightarrow
\mathrm{Bundle}(M/G_0,\prod_k\mathrm{Aut}_G\(X_k\))}i.e.
\ff{\prod_k\mathrm{Vect}_G(M,\rho_k)\longrightarrow
\prod_k\mathrm{Bundle}(M/G_0,\mathrm{Aut}_G\(X_k\))}see \cite[p. 13]{MishQuit}

A bundle $\xi\in\mathrm{Vect}_G(M,\rho)$
is given by transition functions
\ff{\Psi_{\alpha\beta}(x)\in \underset{k}{\prod} \mathrm{Aut}_{G}\(X_k\)}with
the property that there exist
$h_{\alpha,k}(x)\in G_0$ such that
\ff{h^{-1}_{\alpha,k}(x)pr_k\circ\Psi_{\alpha\beta}(x)h_{\alpha,k}(x)}does
not depends on $k$.
Lets show that can be found transition functions
with the property that
\ff{\prod_kpr_k\Psi_{\alpha\beta}(x)=\Delta(a_{\alpha\beta}(x))} for
some cocycle \f{a_{\alpha\beta}(x)}.

Because the group $G_0$ is discrete, for an atlas of connected charts
with connected intersections, we can assume that
\f{pr_k\circ\Psi_{\alpha\beta}(x)=a_{\alpha\beta,k}} and
$h_{\alpha,k}(x)=h_{\alpha,k}\in G_0$ do
not depend on  $x$ and, therefore,
 \ff{h^{-1}_{\alpha,k}(x)pr_k\circ\Psi_{\alpha\beta}(x)h_{\alpha,k}(x)=a_{\alpha\beta}}does
not depend on $k$ nor $x$.
 Let $H_{\alpha,k}\in\mathrm{Aut}_{G}\(X_k\)$ such that
 \f{pr_k(H_{\alpha,k})=h_{\alpha,k}.}Then,
 \ff{\prod_kpr_k(H^{-1}_{\alpha,k}\Psi_{\alpha\beta}(x)H_{\alpha,k})
 =\Delta(a_{\alpha\beta}).} \qed

 \begin{theorem} If the space   $X$ is compact , then
\feq{unionprod}{\mathrm{Bundle}(X,\mathrm{Aut}_G\(X_\rho\))
\approx \bigsqcup_{M\in \mathrm{Bundle}(X,G_0)}
\mathrm{Vect}_G(M,\rho).}
\end{theorem}

\proof We will follow the proof theorem 3 in
\cite[p. 14]{MishQuit}. Given a bundle
$\xi\in \mathrm{Bundle}(X,\mathrm{Aut}_G\(X_\rho\))$
 with transition functions
\ff{\Psi_{\alpha\beta}(x)\in  \mathrm{Aut}_{G}\(X_\rho\)}
we obtain transition functions
\ff{pr\circ i\circ\Psi_{\alpha\beta}(x)\in
\mathrm{Aut}_{G}\(G_0\)\approx G_0,}defining an element
$M\in \mathrm{Bundle}(X,G_0)$ together with
a projection $\xi\mapr{}M$. Changing the fibers
$\mathrm{Aut}_G\(X_\rho\)$ by $X_\rho$,
we obtain an action of the group $G$, that reduces over the base
to the factor group  $G_0$.\qed

Lets rewrite this in terms of homotopy.
\begin{cor} If the space $X$ is compact, then
\feq{unionprod}{[X,\mathrm{Aut}_G\(X_\rho\)]
\approx \bigsqcup_{M\in [X,BG_0]}
\prod_k[M,BU(F_k)].}
\end{cor}\qed

\section{The case when the subgroup $H<G$ is not normal}

Consider an equivariant vector $G$-bundle $\xi$ over the base $M$
\ff{
\begin{array}{c}
          \xi \\
          \mapd{p} \\
          M. \\
        \end{array}
}Let  $H<G$ be a finite subgroup. Assume
that $M$ is the set of fixed points of the conjugation class of this
subgroup, more accurately
\feq{fixconjugate}{M=\underset{[g]\in G/N(H)}{\bigcup} M^{gHg^{-1}},}
and that there is no more fixed points of the conjugation class of $H$ in
the total space of the bundle $\xi$; here we have denoted by
\f{M^H} the set of fixed points of the action of the subgroup
 $H$ over the space $M$, $N(H)$ the normalizer of the group  $H$ in $G$ and
we are using the equality $gM^{H}=M^{gHg^{-1}}$ and the fact that
$lHl^{-1}=gHg^{-1}$
if and only if  $g^{-1}l\in N(H)$.

Lets denote by $\mathfrak{F}_\xi$ the family of subgroups of $G$ having
non-trivial fixed points in the total space of the bundle $\xi$, i.e.
\ff{\mathfrak{F}_\xi=\{K<G| \xi^K\neq \emptyset\}.}
This is a partial ordered set by inclusions and is closed under
the action of the group $G$ by conjugation
\footnote{If $\xi^K\neq \emptyset$, then
                    $\xi^{gKg^{-1}}=g\xi^{K}\neq \emptyset$.}.
Also, the action
\ff{\begin{array}{ccc}G\times\mathfrak{F}_\xi&\mapr{}&\mathfrak{F}_\xi\\
                                              (g,K)                     &\mapsto&gKg^{-1}\\
                                              \end{array}}preserves
the order.

\begin{definition}
We will say that  $H<G$ is the unique, up to
conjugation, maximal subgroup for the $G$-bundle
$\xi$ if every conjugate $gHg^{-1}$ is maximal in
$\mathfrak{F}_\xi$ and there is no more
maximal elements in this family.
\end{definition}

In this section will assume in any case, that $H<G$ is the unique, up to
conjugation, maximal subgroup.

\begin{lemma}\label{voidint}
If $H\neq gHg^{-1}$, then
\ff{M^{H}\cap M^{gHg^{-1}}=\emptyset}
\end{lemma}
\proof\, If there is an
\f{x\in M^{H}\cap M^{gHg^{-1}}}
then, the point \f{x} is fixed under the action of the subgroup generated by
$H$ and $gHg^{-1}$, but this group is not contained in any
of the subgroups of the form $lH^{-1}l,\;l\in G$.\qed

\begin{lemma}\label{disjointbundles}If the condition (\ref{fixconjugate}) holds,
then the $G$-bundle $\xi$ can be presented as a
disjoint union of pair-wise isomorphic bundles with
quasi-free action over the base. More precisely
\ff{\xi=\underset{[g]\in G/N(H)}{\bigsqcup} \xi_{[g]},}where
\ff{\xi_{[g]}=p^{-1}(M^{gHg^{-1}})}
 is a vector bundle with quasi-free
action of the group $N(gHg^{-1})$ and, for every element $g\in G$
  the mapping
\ff{g:\xi^{[1]}\mapr{}g\xi^{[1]}=\xi^{[g]}}
defines an equivariant isomorphism of this bundles, i.e.
the diagram
 \feq{equivariant}{
\begin{array}{ccc}
  N(H)\times \xi_{[1]} &\mapr{} & \xi_{[1]} \\
  \mapd{s_g\times g} &&\mapd{g}\\
  N(gHg^{-1})\times \xi_{[g]}&\mapr{}& \xi_{[g]}\\
\end{array}}commutes, where
\ff{s_g:N(H)\mapr{}N(gHg^{-1})=gN(H)g^{-1},\quad(g,n)\mapsto gng^{-1}.}
\end{lemma}

\proof\,From lemma \ref{voidint} it follows that
 \ff{M=\underset{[g]\in G/N(H)}{\bigsqcup} M^{gHg^{-1}}}and,
 therefore,
 \ff{\xi=\underset{[g]\in G/N(H)}{\bigsqcup} \xi_{[g]}.}Since
the action of  $G$ is fiberwise, we have $g\cdot \xi_{[1]}=\xi_{[g]}$
 for every $g\in G$. Restricting
 the projection $\xi\mapr{} M$ to the space \f{\xi_{[g]}}, we
 obtain the bundle
\ff{\begin{array}{c}
          \xi_{[g]} \\
          \mapd{p} \\
          M^{gHg^{-1}}. \\
          \end{array}}

 The bundle
 $\xi_{[g]}$ has an action of the normalizer $N(gHg^{-1})$:
 \ff{N(gHg^{-1})\times \xi_{[g]}\longrightarrow \xi_{[g]},} i.e.
 \f{\xi_{[g]}} is a $N(gHg^{-1})$-bundle for every $g\in G$.

 Note that group conjugation
$s_g:N(H)\mapr{} N(gHg^{-1})$ defines an isomorphism
 between these groups that fits into the commutative
 diagram
 \ff{
\begin{array}{ccc}
  N(H)\times \xi_{[1]} &\mapr{} & \xi_{[1]} \\
  \mapd{} &&\mapd{}\\
  N(gHg^{-1})\times \xi_{[g]}&\mapr{}& \xi_{[g]}.\\
\end{array}}i.e. $gng^{-1}\cdot gx=g\cdot nx$. This means that
the bundles \f{\xi_{[1]}} and \f{\xi_{[g]}}
are naturally and equivariantly isomorphic.

Evidently, the mappings on the diagram (\ref{equivariant})
do not depend on the elements $n\in N(H)$,
but they depend on the element $g\in G$.

The action of the group
$N(H)$ over the base $M^H$ reduces to the factor group $N(H)/H$:
\ff{
\begin{array}{ccc}
  N(H)\times  \xi_{[1]}&\mapr{} &\xi_{[1]} \\
  \mapd{} &&\mapd{}\\
  N(H)/H\times M^H&\mapr{}& M^H\\
\end{array}}where, considering the maximality of the group $H$,
the action $N(H)/H\times M\mapr{} M$ is free and, by hypothesis,
there is no more fixed of the action of the subgroup
$H$ in the total space of the bundle $\xi$, i.e.
$N(H)$  acts quasi-freely over the base and has normal stationary subgroup $H$.\qed

\begin{definition}
If the condition (\ref{fixconjugate}) holds, we will say that the group
$G$ acts quasi-freely over the bundle $\xi$
with (non-normal) stationary subgroup $H$.
\end{definition}

As we will see in theorem \ref{redtonormal}, for classifying purposes,
it is enough to consider bundles with normal stationary subgroup.

Let $X(\rho)$ be the canonical model for the representation
\f{\rho:H\mapr{} GL(F)} with action of the group $N(H)$. Define
a canonical model  $X(\rho_g)$ for the representation
\ff{\rho_g:gHg^{-1}\mapr{s_{g^{-1}}} H\mapr{\rho} GL(F),}$s_g(n)
=gng^{-1}$. The action of the group $N(gHg^{-1})$ over $X(\rho_g)$ is defined
 using the homomorphism of right $gHg^{-1}$-modules
\ff{u_g:gHg^{-1}\mapr{s_{g^{-1}}}H\mapr{u}N(H)\mapr{s_g} N(gHg^{-1})}by
the formula (\ref{canaction}).

Let
\ff{GX(\rho):=\underset{[g]\in G/N(H)}{\bigsqcup} X(\rho_g)}i.e.
if  $lHl^{-1}=gHg^{-1}$, then the spaces $X(\rho_g)$ and $X(\rho_l)$
coincide.

This notation will be clear after the next lemma.

\begin{lemma}\label{nonnormal}
The group $G$ acts over the space $GX(\rho)$  quasi-freely with
(non-normal) stationary subgroup $H$ and, under this action, the space
$GX(\rho)$ coincides with the orbit of the subspace $X(\rho)$. In particular,
we have the relations
\ff{N(H)\(X(\rho)\)=X(\rho)}and
\ff{(GX(\rho))^{gHg^{-1}}=N(gHg^{-1})/gHg^{-1}.}
\end{lemma}
\proof\,
The action
$G\times GX(\rho)\to GX(\rho)$ is defined
in the following way:
for a fixed $g\in G$ define the mapping
\ff{ g :X(\rho)\mapr{} X(\rho_g) }as
\ff{s_g\times \id :N(H)_0\times F\mapr{}
                       N(gHg^{-1})_0\times F}
(\f{N(H)_0=N(H)/H})
and, if $lHl^{-1}=gHg^{-1}$, then the mapping $l:X(\rho)\mapr{}X(\rho_l)$
is chosen to make the diagramm
\feq{canaction2}{\begin{array}{ccc}
   X(\rho_g)  &\mapr{s_g^{-1}\times \id}& X(\rho)\\
   \|   &                                          & \mapd{l^{-1}g}  \\
  X(\rho_{l}) &\mapr{l^{-1}}                      & X(\rho).
\end{array}}commutative, i.e.
\ff{l=(s_g\times \id)\circ(g^{-1}l)}where
the mapping
\feq{intern}{g^{-1}l:X(\rho)\mapr{}X(\rho)=X(\rho_{g^{-1}l})} is 
the canonical left translation by the element $g^{-1}l\in N(H)$.
\qed

\begin{cor}There is an isomorphism
\feq{isom}{g:\mathrm{Aut}_{N(H)}\(X(\rho)\)\mapr{\approx}
                    \mathrm{Aut}_{N(gHg^{-1})}\(X(\rho_g)\)}that depends
only on the class $[g]\in G/N(H)$.
\end{cor}

\proof\, We have a diagram  (\ref{equivariant}) for $\xi=GX(\rho)$.
Such a diagram always induces an isomorphism
\ff{\mathrm{Aut}_{N(H)}\(X(\rho)\)\mapr{\approx}
                    \mathrm{Aut}_{N(gHg^{-1})}\(X(\rho_g)\)}by
the rule
\ff{\mathbf{A}\mapsto g\mathbf{A}g^{-1}}and,
if $l\in [g]\in  G/N(H)$ then $l^{-1}g\in N(H)$ commutes with
$\mathbf{A}\in\mathrm{Aut}_{N(H)}\(X(\rho)\)$. Therefore
\ff{g\mathbf{A}g^{-1}=g(g^{-1}l)(l^{-1}g)\mathbf{A}g^{-1}=g(g^{-1}l)\mathbf{A}(l^{-1}g)g^{-1}=l\mathbf{A}l^{-1}.} \qed

\begin{definition}
The space $GX(\rho)$ is called the canonical model for the case
when the subgroup  $H<G$ is not normal.
\end{definition}

\begin{lemma}
\feq{isoaut}{\mathrm{Aut}_{G}\(GX(\rho)\)\approx
          \mathrm{Aut}_{N(H)}\(X(\rho)\)}
\end{lemma}

\proof\,
By definition, an element of the group \f{\mathrm{Aut}_G\(X\)}
is an equivariant mapping $\mathbf{A}^{a}$ such that the pair
$(\mathbf{A}^{a},a)$
defines the commutative diagram
\ff{
\begin{array}{ccc}
   X &\mapr{\mathbf{A}^{a}}& X\\
  \mapd{} &&\mapd{}\\
   G/H&\mapr{a}&G/H, \\
\end{array}
}that
commutes with the canonical action, i.e. the mapping
 \f{a\in \mathrm{Aut}_G(G/H)} satisfies the condition
\ff{a\in \mathrm{Aut}_G(G/H)\approx N(H)/H,\quad a[g]=[ga],\;[g]\in N(H)/H.}
Therefore,
$\mathbf{A}^{a}=(A^{a}[g])_{[g]\in N(H)/H}\in \mathrm{Aut}_{N(H)}(X(\rho))$.

The value of the operators
$(A^{a}[g])_{[g]\in G/H}$
can be calculated in terms of the operator
$A^{a}[1]$ as in lemma 2 from \cite[p. 9]{MishQuit}.
\qed

Denote by
$\widetilde{\mathrm{Vect}}_G(M,\rho)$ the category of vector bundles
with quasi-free action of the group $G$ over the base $M$.

\begin{theorem}\label{redtonormal}
$\widetilde{\mathrm{Vect}}_G(M,\rho)\approx \mathrm{Vect}_{N(H)}(M^H,\rho)$.
\end{theorem}

\proof\,
From lemma \ref{disjointbundles} follows that the bundles
$\xi_{[1]}$ and $\xi_{[g]}$ equivariantly isomorphic and
are given by mappings
\ff{
M^{gHg^{-1}}/N(gHg^{-1})_0\mapr{} B\mathrm{Aut}_{N(gHg^{-1})}\(X(\rho_g)\),
}and
\ff{
M^H/N(H)_0\mapr{} B\mathrm{Aut}_{N(H)}\(X(\rho)\),
}that can be put in the
commutative diagram
\ff{\begin{array}{ccc}
   M^{H}/N(H)_0&\mapr{}& B\mathrm{Aut}_{N(H)}\(X(\rho)\)\\
   \mapd{\bar g}& &\mapd{\bar g}\\
  M^{gHg^{-1}}/N(gHg^{-1})_0&\mapr{}& B\mathrm{Aut}_{N(gHg^{-1})}\(X(\rho_g)\).
\end{array}
}Here, $g:\xi^{H}\mapr{}\xi^{gHg^{-1}}$ is the action over the bundle
 $\xi$. The arrow on the right side is induced by the isomorphism
(\ref{isom}) and does not depend on the element $g\in [g]\in G/N(H)$.
\qed


\begin{thebibliography}{20}
\bibitem{MishQuit}     Mishchenko A.S., Morales Mel\'endez, Quitzeh.
                                     \textit{Description of  $G$-bundles over
                                     $G$-spaces with quasi-free proper action
                                      of discrete group}
\bibitem{Mishchenko} Luke G., Mishchenko A. S., \textit{Vector Bundles And Their Applications.}
                                       Kluwer Academic Publishers Group (Netherlands), 1998.
                                       ISBN: 9780792351542
\bibitem{Connor} P. Conner, E. Floyd. \textit{Differentiable periodic maps.} Berlin, Springer-Verlag 1964.
\bibitem{Palais} Palais R.S. \textit{On the Existence of Slices for Actions of Non-Compact Lie Groups}
                            Ann.  Math., 2nd Ser., Vol. 73, No. 2. (1961), pp. 295-323.
\bibitem{Atiyah} Atiyah M.F., \textit{K-theory.} Benjamin, New York, (1967).
\bibitem{Serre}Serre J.P., \textit{Representations line\'{a}ires des groupes finis.} Hermann, Paris. 1967.
\bibitem{Levine} Levine M.,  Serp\'{e} C.,\textit{On a spectral sequence for equivariant K-theory}
                                                                        K-Theory (2008) 38 pp. 177–222
\bibitem{BeylTappe} Beyl F.R.,  Tappe J.,\textit{Groups Extensions, Representations, and the Schur Multiplicator.}
                                                                        Springer-Verlag (Berlin Heidelberg), 1982.
                                                                        ISBN 354011954X
\bibitem{Brown}  Brown K.S. \textit{Cohomology of groups.}
                                                   Springer-Verlag (New York Heidelberg Berlin), 1982.
\bibitem{EilenbergMacLane}  Eilenberg S., MacLane S.
                                                    \textit{Cohomology theory in abstract groups. II. Group extensions with a
                                                      non-abelian kernel.}
                                                      Ann. Math., 1947, (2) 48, p. 326—341.

\end{thebibliography}
\end{document}